\documentclass[12pt,a4paper,leqno]{article}
\usepackage{etoolbox,tocloft} \usepackage{comment}
\usepackage{amsfonts}
\usepackage[a4paper,includeheadfoot,margin=0.87in]{geometry}
\usepackage{times}
\usepackage{amsthm,amssymb,mathtools}
\usepackage{amsmath,xypic}
\usepackage[export]{adjustbox} 
\usepackage[cal=boondoxo,calscaled=.96,scr=rsfs]{mathalpha} \usepackage{multicol}
\usepackage{hyphenat}
\usepackage[usenames,dvipsnames]{color}
\usepackage{rotating}
\usepackage{tikz-cd}
\usepackage{longtable}
\usepackage{caption}
\usepackage[most]{tcolorbox}

\usepackage{epigraph}
\usepackage{enumitem}
\setlist[enumerate]{listparindent=0.5in}
\usepackage{mathrsfs}
\newcommand{\be}{\begin{equation}}
\newcommand{\ee}{\end{equation}}
\newcommand{\bes}{\begin{equation*}}
\newcommand{\ees}{\end{equation*}}
\newcommand{\bea}{\begin{eqnarray}}
\newcommand{\eea}{\end{eqnarray}}
\newcommand{\beas}{\begin{eqnarray}}
\newcommand{\eeas}{\end{eqnarray}}
\newcommand{\ben}{\begin{note}}
\newcommand{\een}{\end{note}}
\newcommand{\bexl}{\vskip0.1em\noindent\hrulefill\vskip1em\begin{ExerciseList}}
\newcommand{\eexl}{\end{ExerciseList}\hrulefill}

\newcommand{\bthm}{\begin{theorem}}
\newcommand{\ethm}{\end{theorem}}
\newcommand{\bpro}{\begin{prop}}
\newcommand{\epro}{\end{prop}}
\newcommand{\bcor}{\begin{corollary}}
\newcommand{\ecor}{\end{corollary}}
\newcommand{\bcon}{\begin{conjecture}}
\newcommand{\econ}{\end{conjecture}}
\newcommand{\bp}{\begin{proof}}
\newcommand{\ep}{\end{proof}}
\newcommand{\blem}{\begin{lemma}}
\newcommand{\elem}{\end{lemma}}
\newcommand{\bn}{\begin{note}}
\newcommand{\en}{\end{note}}
\newcommand{\benum}{\begin{enumerate}}
\newcommand{\eenum}{\end{enumerate}}
\newcommand{\bed}{\begin{defn}}
\newcommand{\eed}{\end{defn}}
\newcommand{\brem}{\begin{remark}}
\newcommand{\erem}{\end{remark}}

\newcommand{\btik}{\begin{tikzpicture}\begin{axis}[scale=0.5,axis y line=center, axis x line=middle]}
\newcommand{\etik}{\end{axis}\end{tikzpicture}}

\newcommand{\upperRomannumeral}[1]{\uppercase\expandafter{\romannumeral#1}}

 \usepackage{graphicx}
\usepackage[pagewise]{lineno}
\usepackage{stackengine}
\usepackage{stmaryrd}
\usepackage[normalem]{ulem}

\usepackage{titlesec}
\usepackage{url}	
\usepackage{natbib}
\let\cite=\citep
\usepackage{fancyhdr}

\usepackage[colorlinks,citecolor=blue,colorlinks=true,hyperindex, citecolor=blue, urlcolor=blue]{hyperref}
\usepackage[capitalise,noabbrev]{cleveref}

\newtoggle{draft}
\togglefalse{draft}

\newtheoremstyle{mytheorem}
{}{}{\itshape}{}{\bfseries\sffamily}{}{\newline}{}

\theoremstyle{mytheorem}
\newtheorem{theorem}[equation]{Theorem}      \newtheorem{lemma}[equation]{Lemma}          \newtheorem{corollary}[equation]{Corollary}  \newtheorem{proposition}[equation]{Proposition}

\newtheorem{conj}[equation]{Conjecture}

\newtheorem{defn}[equation]{Definition}

\newtheorem{remark}[equation]{Remark}

\numberwithin{equation}{section}

\let\oldproofname=\proofname
\renewcommand{\proofname}{{\bfseries\sffamily\textup{\oldproofname}}}

\renewcommand{\thesubsection}{\S~\arabic{section}.\arabic{subsection}}

\titleformat{\subsection}[runin]{\normalfont\bfseries}{\thesubsection}{.5em}{}[{\ \ }]
\titlespacing{\subsection}{0pt}{1.5ex plus .1ex minus .2ex}{0pt}
\titleformat{\subsubsection}[runin]{\normalfont\bfseries}{\thesubsubsection}{.5em}{}[{\ \ }]
\titlespacing{\subsubsection}{0pt}{1.5ex plus .1ex minus .2ex}{0pt}

\crefname{section}{§}{§§}
\crefname{subsection}{§}{§§}
\crefname{subsubsection}{§}{§§}

\usepackage{microtype}
\usepackage[margin=15pt,font=small,labelfont={bf,sf}]{caption}

\newcommand{\A}{\mathscr{A}}

\newcommand{\bF}{{\bar{F}}}
\newcommand{\bQ}{{\bar{\Q}}}

\newcommand{\C}{{\mathbb C}}

\newcommand{\F}{{\mathbb F}}

\newcommand{\N}{\mathscr{N}}

\newcommand{\Q}{{\mathbb Q}}

\renewcommand{\O}{{\mathscr O}}
\renewcommand{\P}{{\mathbb P}}
\renewcommand{\wp}{{\mathfrak p}}

\renewcommand{\bpro}{\begin{proposition}}
	\renewcommand{\epro}{\end{proposition}}
\renewcommand{\bcon}{\begin{conj}}
	\renewcommand{\econ}{\end{conj}}

\title{Final Report on the Mochizuki-Scholze-Stix Controversy}
\date{March 12, 2025}
\author{Kirti Joshi}

\newcommand{\Address}{\bigskip\noindent{\footnotesize\textsc{{Math. department, University of Arizona, 617 N Santa Rita, Tucson
		85721-0089, USA.}}\par\nopagebreak 
}}

\begin{document}
	\maketitle
\epigraphwidth0.65\textwidth

\lhead{}

\iftoggle{draft}{\pagewiselinenumbers}{\relax}
\newcommand{\act}{\curvearrowright}
\newcommand{\lmp}{{\Pi\act\Ot}}
\newcommand{\lmpi}{{\lmp}_{\int}}
\newcommand{\lmpf}{\lmp_F}
\newcommand{\Om}{\O^{\times\mu}}
\newcommand{\Omf}{\O^{\times\mu}_{\bF}}
\renewcommand{\N}{\mathbb{N}}
\newcommand{\yoga}{Yoga}
\newcommand{\gl}[1]{{\rm GL}(#1)}
\newcommand{\bK}{\overline{K}}
\newcommand{\reptrip}{\rho:G_K\to\gl V}
\newcommand{\reptripp}[1]{\rho\circ\alpha:G_{\ifstrempty{#1}{K}{{#1}}}\to\gl V}
\newcommand{\benumlab}{\begin{enumerate}[label={{\bf(\arabic{*})}}]}
\newcommand{\ord}{\mathop{\rm ord}\nolimits}	
\newcommand{\kcs}{K^\circledast}
\newcommand{\lcs}{L^\circledast}
\renewcommand{\A}{\mathbb{A}}
\newcommand{\bfq}{\bar{\mathbb{F}}_q}
\newcommand{\tripod}{\P^1-\{0,1728,\infty\}}

\newcommand{\vseq}[2]{{#1}_1,\ldots,{#1}_{#2}}
\newcommand{\anab}[4]{\left({#1},\{#3 \}\right)\anabelmap\left({#2},\{#4 \}\right)}

\newcommand{\gln}{{\rm GL}_n}
\newcommand{\glo}[1]{{\rm GL}_1(#1)}
\newcommand{\glt}[1]{{\rm GL_2}(#1)}

\newcommand{\iut}{\cite{mochizuki-iut1, mochizuki-iut2, mochizuki-iut3,mochizuki-iut4}}
\newcommand{\topics}{\cite{mochizuki-topics1,mochizuki-topics2,mochizuki-topics3}}

\newcommand{\linv}{\mathfrak{L}}
\newcommand{\bedef}{\begin{defn}}
\newcommand{\eedef}{\end{defn}}
\renewcommand{\act}[1][]{\overset{#1}{\curvearrowright}}
\newcommand{\bfx}{\overline{F(X)}}
\newcommand{\anabelmap}{\leftrightsquigarrow}
\newcommand{\ban}[1][G]{\mathscr{B}({#1})}
\newcommand{\pit}{\Pi^{temp}}
 
 \newcommand{\bL}{\overline{L}}
 \newcommand{\bkm}{\bK_M}
 \newcommand{\vbk}{v_{\bK}}
 \newcommand{\vbkm}{v_{\bkm}}
\newcommand{\ocs}{\O^\circledast}
\newcommand{\ot}{\O^\triangleright}
\newcommand{\ocsk}{\ocs_K}
\newcommand{\otk}{\ot_K}
\newcommand{\ok}{\O_K}
\newcommand{\oko}{\O_K^1}
\newcommand{\oks}{\ok^*}
\newcommand{\Qpb}{\overline{\Q}_p}
\newcommand{\Qpbh}{\widehat{\overline{\Q}}_p}
\newcommand{\tr}{\triangleright}
\newcommand{\ocpt}{\O_{\C_p}^\tr}
\newcommand{\ocpf}{\O_{\C_p}^\flat}
\newcommand{\sG}{\mathscr{G}}
\newcommand{\sX}{\mathscr{X}}
\newcommand{\sxfe}{\sX_{F,E}}
\newcommand{\sxfep}{\sX_{F,E'}}
\newcommand{\sxcpte}{\sX_{\cpt,E}}
\newcommand{\sxcptep}{\sX_{\cpt,E'}}
\newcommand{\loglt}{\log_{\sG}}
\newcommand{\fc}{\mathfrak{t}}
\newcommand{\ku}{K_u}
\newcommand{\kup}{\ku'}
\newcommand{\kt}{\tilde{K}}
\newcommand{\sGpf}{\sG(\O_K)^{pf}}
\newcommand{\hgm}{\widehat{\mathbb{G}}_m}
\newcommand{\bE}{\overline{E}}
\newcommand{\sY}{\mathscr{Y}}
\newcommand{\syfe}{\mathscr{Y}_{F,E}}
\newcommand{\syfep}{\mathscr{Y}_{F,E'}}
\newcommand{\syfqp}[1]{\mathscr{Y}_{\cptl{#1},\Q_p}}
\newcommand{\syfqpe}[1]{\mathscr{Y}_{{#1},E}}
\newcommand{\syfqpep}[1]{\mathscr{Y}_{{#1},E'}}
\newcommand{\fJ}{\mathfrak{J}}
\newcommand{\fM}{\mathfrak{M}}
\newcommand{\locvar}{local arithmetic-geometric anabelian variation of fundamental group of $X/E$ at $\wp$}
\newcommand{\fjxep}{\fJ(X,E,\wp)}
\newcommand{\fjxe}{\fJ(X,E)}
\newcommand{\fpc}[1]{\widehat{{\overline{\F_p(({#1}))}}}}
\newcommand{\cpt}{\C_p^\flat}
\newcommand{\cpti}{\C_{p_i}^\flat}
\newcommand{\cptl}[1]{\C_{p,{#1}}^\flat}
\newcommand{\fja}[1]{\fJ^{\rm arith}({#1})}
\newcommand{\ainfe}{A_{\inf,E}(\O_F)}
\renewcommand{\ainfe}{W_{\O_E}(\O_F)}
\newcommand{\gmh}{\widehat{\mathbb{G}}_m}
\newcommand{\sE}{\mathscr{E}}
\newcommand{\bpi}{B^{\varphi=\pi}}
\newcommand{\bpip}{B^{\varphi=p}}
\newcommand{\onto}{\twoheadrightarrow}

\newcommand{\cpmax}{{\C_p^{\rm max}}}
\newcommand{\xan}{X^{an}}
\newcommand{\yan}{Y^{an}}
\newcommand{\bPi}{\overline{\Pi}}
\newcommand{\bPit}{\bPi^{\rm{\scriptscriptstyle temp}}}
\newcommand{\Pit}{\Pi^{\rm{\scriptscriptstyle temp}}}
\renewcommand{\pit}[1]{\Pi^{\scriptscriptstyle temp}_{#1}}
\newcommand{\pitk}[2]{\Pi^{\scriptscriptstyle temp}_{#1;#2}}
\newcommand{\pio}[1]{\pi_1({#1})}
\newcommand{\fTeich}{\widetilde{\fJ(X/L)}}
\newcommand{\ssep}{\S\,} \newcommand{\vphi}{\varphi}
\newcommand{\sgt}{\widetilde{\sG}}
\newcommand{\sxqp}{\mathscr{X}_{\cpt,\Q_p}}
\newcommand{\syQp}{\mathscr{Y}_{\cpt,\Q_p}}

\setcounter{tocdepth}{2}
\tableofcontents

\newcommand{\mywork}[1]{\textcolor{teal}{#1}}

\togglefalse{draft}
\newcommand{\FF}{\cite{fargues-fontaine}}
\iftoggle{draft}{\pagewiselinenumbers}{\relax}

\newcommand{\attportion}{Sections~\cref{se:number-field-case}, \cref{se:construct-att}, \cref{se:relation-to-iut}, \cref{se:self-similarity} and \cref{se:applications-elliptic}}

\newcommand{\Pib}{\overline{\Pi}}
\newcommand{\four}{Sections~\cref{se:grothendieck-conj}, \cref{se:untilts-of-Pi}, and \cref{se:riemann-surfaces}}

\numberwithin{equation}{subsection}
\newcommand{\omu}{\O_{\bQ_p}^{\mu}}
\newcommand{\lmod}{L_{\rm mod}}
\newcommand{\ttheta}{\widetilde{\Theta}}
\newcommand{\tthetaj}[1]{\ttheta_{Joshi,#1}}
\newcommand{\tthetam}[1]{\ttheta_{Mochizuki,#1}}

\newcommand{\moccor}{\cite[Corollary 3.12]{mochizuki-iut3}}
\newcommand{\thetam}{{\ttheta}_{Mochizuki}}
\newcommand{\thetamp}{{\ttheta}_{Mochizuki,p}}
\newcommand{\thetaj}{{{\ttheta}_{Joshi}}}
\newcommand{\thetajp}{{{\ttheta}_{Joshi,p}}}
\newcommand{\sM}{\mathscr{M}}
\newcommand{\pib}{\overline{\Pi}}
\newcommand{\bN}{\mathbb{N}}
\newcommand{\sD}{\mathscr{D}}
\newcommand{\sF}{\mathscr{F}}
\newcommand{\sL}{\mathscr{L}}

\newcommand{\bdrp}{B_{dR}^+}
\newcommand{\bdr}{B_{dR}}

\newcommand{\iutthr}{\cite{mochizuki-iut1,mochizuki-iut2, mochizuki-iut3}}

\newcommand{\thetajpp}{{\widehat{\Theta}}_{Joshi,p}}
\newcommand{\thetaja}{{\widehat{\Theta}}_{Joshi}}
\newcommand{\thetajpph}{{\widehat{\widehat{\Theta}}}_{Joshi,p}}
\newcommand{\thetajppa}{{\widehat{\widehat{\Theta}}}_{Joshi}}

\newcommand{\ells}{{\ell^*}}

\definecolor{darkmidnightblue}{rgb}{0.0, 0.2, 0.4}
\definecolor{carmine}{rgb}{0.59, 0.0, 0.09}
\newtcbox{\mybox}[1][red]{on line,
	arc=0pt,outer arc=0pt,colback=white,colframe=darkmidnightblue,
	boxsep=0pt,left=0pt,right=0pt,top=2pt,bottom=2pt,
	boxrule=0pt,leftrule=1pt, rightrule=1pt,bottomrule=1pt,toprule=1pt}
\newcommand{\tmb}[1]{\mybox{#1}}
\newcommand{\present}{the series of papers (\cite{joshi-teich,joshi-untilts, joshi-teich-estimates,joshi-teich-def,joshi-teich-rosetta,joshi-teich-abc-conj} and \cite{joshi-anabelomorphy,joshi-formal-groups})}

\newtheorem{note}{Note}
\numberwithin{note}{subsubsection}

\newtcolorbox[auto counter,crefname={rosetta Stone Fragment}{rosetta Stone Fragements}]{rosetta}[2][]{colback=white,coltitle=black,colframe=white!25!brown,fonttitle=\bfseries, 
	title=Table \thetcbcounter\ #2,#1}
\nocite{joshi-formal-groups,joshi-anabelomorphy,joshi-gconj,joshi-untilts-2020,joshi-untilts,joshi-teich,joshi-teich-summary-comments,joshi-teich-estimates,joshi-teich-quest,joshi-teich-def,joshi-teich-rosetta,joshi-teich-abc-conj}
\addtocounter{section}{1} 
\renewcommand{\footnotesize}{\fontsize{10pt}{12pt}\selectfont}
\setlength{\footnotesep}{\baselineskip}
\subsection{The purpose of this document}
This brief document records my mathematical findings vis-a-vis \iut\ and \cite{scholze-stix,scholze-review}. My findings are tabulated in \ref{ss:tabulation-of-findings} (Tables \ref{tab:rosetta-stone1}, \ref{tab:rosetta-stone2}). As Table \ref{tab:rosetta-stone2} shows, every assertion of \cite{scholze-stix} and \cite{scholze-review} is mathematically false. On the other hand, Mochizuki's proof is also incomplete (see \ref{se:why}). A robust version of the theory claimed by Mochizuki is provided by my work. Proofs of all my mathematical assertions in this context may be found in the revised versions\footnote{This version of [Constr. IV] was sent to Mochizuki and Scholze in Nov 2024.} (Feb 2025) of \present.

My previous report on this topic was released in June 2024 \cite{joshi-report} which includes a timeline of events (reading \cite{joshi-report} is still recommended). The conclusion of this report (\ref{ss:conclusion}) replaces the conclusion of that earlier report. 
\subsection{Why is Mochizuki's proof incomplete?}\label{se:why}
Mochizuki's remarkable claim of the proof of the $abc$-conjecture rests on an  astounding assertion that there exists a Teichm\"uller Theory\footnote{Teichm\"uller Theory should not be confused with Riemann's Moduli Theory of Riemann surfaces.} of number fields 
(that is why the phrase Teichm\"uller Theory appears in the title of \iut). Mochizuki correctly surmised that such a theory exists, however, that is not enough to prove the existence of such a theory nor are Mochizuki's anabelian geometry methods adequate for  demonstrating its existence. 

At the very center of the issue is that Mochizuki's quantification of what it means to be an \textit{Arithmetic Holomorphic Structure} is mathematically inadequate to quantitatively assert that one has two or more such structures. This last point is needed because in the proof of the $abc$-conjecture one compares and averages over arithmetic quantities arising from many such structures, and so it is important to unambiguously prove the existence of many such structures. In particular, without clearly proving this, there is no way to provide a complete proof of \cite[Theorem 3.11]{mochizuki-iut3} (let alone \cite[Corollary 3.12]{mochizuki-iut3}). 

Peter Scholze and Jakob Stix recognized this problem (2018) -- but they extrapolated and argued (incorrectly)  that many such structures cannot exist (see \ref{ss:tabulation-of-findings}). 

My work  (\cite{joshi-teich}, \cite{joshi-untilts}) provides a precise definition of `Arithmetic Holomorphic Structures' and this  allows me to 
\benumlab
\item show that these arise from rigid analytic spaces i.e. arise from $p$-adic holomorphic (or analytic) functions--hence one has the $p$-adic analog of classical Teichm\"uller Theory;
\item explicitly exhibit the existence of many such structures and, 
\item allows one to elaborate and transparently prove the properties (of such structures) which Mochizuki has claimed in his work and requires in his proofs.  
\eenum
I have laid out my theory in the `Construction of Arithmetic Teichm\"uller Spaces' series of papers. 

There is one important point which needs to be clearly understood: Mochizuki has argued that his proof exists because of subtle aspects of Anabelian Geometry (and group theory surrounding fundamental groups). My finding is that this is mathematically not the case.  My finding is that the theory exists for a subtler and deeper reason:

 \textit{Arithmetic, both local and global, is far richer and occurs in  many topologically distinct avatars than has been previously imagined  (\cite{joshi-teich-def}).} 

Mochizuki, no doubt, surmised this standing atop \cite[Theorem 1.9]{mochizuki-topics3} (and proceeded to build his theory). But the main lacuna in the said theorem is that it is quantitatively ineffective in providing \iut\ with a way of exhibiting many distinct arithmetic holomorphic structures or establishing the properties of these structures. 

My discovery (\cite{joshi-untilts}, \cite{joshi-teich}, \cite{joshi-teich-def}) is that the Theory of Arithmetic Teichm\"uller Spaces, which includes Mochizuki's Inter-Universal Teichm\"uller Theory as a special case, is birthed by the vast and profound richness of $p$-adic Arithmetic itself. 

That there are many inequivalent versions of $p$-adic arithmetic has been staring at us since \cite{schmidt33}, \cite{kaplansky42}, \cite{matignon84}, \cite{kedlaya18}. This manifests itself in the category of algebraically closed perfectoid fields of characteristic zero and prime residue characteristics\footnote{Note that one must work with  untilts--it is not enough to work with algebraically closed perfectoid fields.} (\cite{scholze12-perfectoid-ihes}, \cite{fargues-fontaine}) and its consequences for Mochizuki's Theory appear in my works beginning with \cite{joshi-untilts-2020}, \cite{joshi-teich}, \cite{joshi-untilts} and  \cite{joshi-teich-def}, where this leads to the global version, which precisely says that a fixed number field is topologically deformable and that there are many such deformations of global arithmetic (as asserted by Mochizuki but not proved by him). The category of algebraically closed perfectoid fields with prime residue characteristics (there is one for every prime number $p$ and all of these are simultaneously needed) is also the (only) source of (arbitrary) geometric base-points for tempered fundamental groups, required by \textit{Mochizuki's Key Principle of Inter-Universality} \cite[\ssep I3, Pages 25-26]{mochizuki-iut1} which lies at the foundation of Mochizuki's Theory. 

But these categories (of algebraically closed perfectoid fields of characteristic zero and residue characteristic $p>0$)  do not exist because of Anabelian Geometry. \textit{They simply exist}. 

That is why, it would be completely incorrect to say, as many around Mochizuki have repeatedly said for the past decade, that \iut\ is about Anabelian Geometry. [Mochizuki  recognized this issue (post publication of his papers) and tried to supplement his proofs by means of the `change of logic' argument of \cite{mochizuki-essential-logic}. But that solution is not satisfactory and some of its arguments, notably analogy with the construction of projective spaces, is not mathematically useful because one wants to average over distinct arithmetic holomorphic structures in \cite[Theorem 3.11, Corollary 3.12]{mochizuki-iut3}]. 

More importantly, this argument of \cite{mochizuki-essential-logic} is not only flawed (because it simply declares the existence of distinct arithmetic holomorphic structures), but mathematically superfluous. There is a canonical definition of arithmetic holomorphic structures (\cite{joshi-teich}, \cite{joshi-untilts}) and avoiding it in hopes of avoiding algebraically closed perfectoid fields (and the relevant mathematics) serves no mathematical purpose at this juncture (and contradicts Mochizuki's Key Principle of Inter-Universality). In particular, there is one and only one theory  which can be constructed using Mochizuki's Key Principle and that theory is described in my papers. [\cite{joshi-teich-rosetta} provides a `Rosetta Stone' to facilitate a parallel reading of the two theories.] 

The idea that there exist many topologically inequivalent versions of arithmetic is truly remarkable and an extremely subtle one (frankly, most mathematicians who engaged with Mochizuki's proof have missed it completely)  and will take some time to sink in. It furthers the Dedekind-Weber analogy between Riemann surfaces and number fields by establishing a Teichm\"uller Theory of number fields.  At this juncture, there is little room for any debate about its existence (\cite{joshi-teich-def}) and this existence underpins Mochizuki's proof of the $abc$-conjecture. Hence its precise quantification is absolutely essential (for validating Mochizuki's claim). I am not prescient enough to see all its consequences. But \cite{joshi-teich-def} shows (going beyond Mochizuki's original claim) that one even has a topological space of deformations of a fixed number field and a Frobenius morphism on this space. This opens up the possibility of treating problems in Diophantine geometry from a point of view of Arithmetic Dynamics using this space and its Frobenius morphism--to some extent, this Arithmetic Dynamics aspect is already at work in the proof of the $abc$-conjecture  \cite{joshi-teich-abc-conj} and \cite{mochizuki-iut4}.

\subsection{Tabulation of Findings}\label{ss:tabulation-of-findings}
The following tables (Tables \ref{tab:rosetta-stone1} and \ref{tab:rosetta-stone2}) provide a comparison of hypotheses of Mochizuki and Joshi and the omissions of these hypothesis by Scholze-Stix leading to incorrect mathematical conclusions.
\begin{rosetta}[label=tab:rosetta-stone1,grow to right by=1cm,grow to left by=1.5cm]{{Comparison of claims of Mochizuki, Scholze-Stix and Joshi}}
\centering{\small
		\begin{tabular}{|p{0.9in}|p{1.9in}|p{1.5in}|p{2in}|}\hline 
			object	&	Mochizuki	& Scholze-Stix & Joshi \\ 
			\hline 
			geom. base point &	Central role of arbitrary geometric base-points as a proxy for deformations of arithmetic \cite[\ssep I3, Page 25]{mochizuki-iut1}	& ignored & included (in the data of an arithmeticoid) as arbitrary alg. closed perfectoid fields \cite{joshi-teich,joshi-untilts} \\ 
			\hline 
			how is this used?	&	domain and codomain of all key operations refer to distinct geom. base-points \cite[\ssep I3, Page 25]{mochizuki-iut1}	& incorrectly identify the domains and codomains leading to incorrect conclusions & naturally show Mochizuki's requirements \cite{joshi-teich-rosetta} \\ 
			\hline 
			Valuation data & encoded in realified Frobenioids \cite{mochizuki-iut1} & no mention of distinct valuation data & works with valuation data instead of Frobenioids \cite{joshi-teich,joshi-untilts,joshi-teich-def}\\
			\hline
			Why needed? & for computing local and global arith. degrees & ignored & as in Mochizuki but without using Frobenioids \cite{joshi-teich-def,joshi-untilts}\\
			\hline
			log-Links & log-Links aka Mochizuki's proxy for Frobenius (at each prime) \cite{mochizuki-iut3} & ignored & works with the Frobenius morphism instead of a proxy \cite{joshi-teich-estimates,joshi-teich-def}.\\
			\hline
			Theta-Links & Central role of Theta-Links Claimed via theory of Frobenioids \cite{mochizuki-iut1} & argue that such objects cannot exist & demonstrates the existence both at local and global level. \cite{joshi-teich-estimates,joshi-teich-def,joshi-teich-rosetta}\\
			\hline
			distinct Arith. Holomorphic Structures & Asserted in \cite{mochizuki-iut1} but existence is not clearly established & declare that these cannot exist & Demonstrates the existence and deformation property via arbitrary alg. closed perfectoid fields or equivalently by using arbitrary geom. base-points \cite{joshi-teich,joshi-untilts,joshi-teich-def}\\
			\hline
		\end{tabular} 
	}
	\tcblower
	{\small
		\noindent\textcolor{red}{{\bf NOTE}}
		\benumlab
		\item Despite the fact that Mochizuki asserts his Key Principle of Inter-Universality \cite[\ssep I3, Page 25]{mochizuki-iut1}, there is no mention of the required input data in \iut, and on two separate occasions, Mochizuki denied the relevance of alg. closed perfectoid fields to \iut.
		\item One cannot build a valid theory, as \iut\ claims to do, by requiring arbitrary geometric base-points for tempered fundamental groups but work solely over $\bQ_p$. This is the reason for my mathematical objections to \iut. My work fixes this central issue and builds a robust (and even more general) theory of Arithmetic Teichm\"uller Spaces, which strictly adheres to Mochizuki's Key Principle of Inter-Universality \cite[\ssep I3, Page 25]{mochizuki-iut1}, and with all properties claimed by Mochizuki for his IUTT.
		\eenum
	}
\end{rosetta}

\newpage
\begin{rosetta}[label=tab:rosetta-stone2,grow to right by=1cm,grow to left by=1.5cm]{{Status of various assertions of \cite{scholze-stix}, \cite{scholze-review}}}
\centering{\small
		\begin{tabular}{|p{2in}|p{2.9in}|p{1.5in}|}\hline 
			Claims of section in \cite{scholze-stix}	& Disproved by Joshi in the \present & Mochizuki (for comparison) \\ 
			\hline 
			Sect. 2.1.2 (and also \cite{scholze-review})
			asserts that distinct
			\benumlab 
			\item Hodge Theaters
			\item \'Etale Pictures
			\item Frobenius Pictures
			\eenum 
			do not exist & \benumlab 
			\item distinct Hodge-Theaters exist \cite[Theorem 10.11.5.1]{joshi-teich-rosetta};
			\item  distinct \'etale Pictures exist \cite[Proposition 8.3.1.1]{joshi-teich-rosetta}; 
			\item distinct Frobenius Pictures exist \cite[Proposition 8.3.1.2]{joshi-teich-rosetta}
			\eenum	& asserts existence of all three (without proof) \\ 
			\hline 
			Remark 9	&  existence of isomorphs \cite[Theorem 2.9.1]{joshi-teich}, \cite[Theorem 4.8]{joshi-untilts} and arithmetic holomorphic structures \cite[Definition 5.1]{joshi-untilts}	 & asserts existence (without proof)\\ 
			\hline 
			Section 2.1.4 & False by \cite[Theorem 5.10.1]{joshi-teich-def} & is tacitly asserted (without proof) \\
			\hline
			Section 2.1.5	& distinct prime-strips exist \cite[\ssep 8.7.1 and Theorem 8.8.3]{joshi-teich-rosetta} & asserts the existence (without proof)  \\
			\hline
			Section 2.1.6 (argues $\mathfrak{log}$-Link must be irrelevant)	& $\mathfrak{log}$-Link = global Frobenius and it is constructed in \cite[\ssep 5]{joshi-teich-def}  (and is non-trivial), $\mathfrak{log}$-Link aspect is detailed in \cite[\ssep 8.9]{joshi-teich-rosetta}  & asserts properties in \cite{mochizuki-topics3} (without proof) \\
			\hline
			Section 2.1.7, 2.1.8		 & Each holomorphoid provides $q$ and $\Theta$-Pilot object \cite[Remark 6.4.2.2 and \ssep 6.10.1]{joshi-teich-rosetta}   & asserts existence (without proof)\\
			\hline
			Section 2.1.9 ($\Theta_{gau}$-Link cannot be non-trivial)		 & Construction and properties \cite[Theorem 4.2.2.1 and \ssep 4.2.3 and \ssep 8.10 ]{joshi-teich-rosetta}  & asserts non-triviality (without proof) \\
			\hline
			Section 2.2 (main disproof argument)	 & This argument of \cite{scholze-stix} fails because the $\Theta_{gau}$-Link does come with a non-trivial $j^2$-scaling factor for $j=1,\ldots,\frac{\ell-1}{2}$ and also is equipped with a non-trivial Galois action  \cite[Theorem 4.2.2.1 and \ssep 4.2.3 and \ssep 8.10 ]{joshi-teich-rosetta}  &  \\
			\hline
		\end{tabular} 
	}
	\tcblower
	{\small
		\noindent\textcolor{red}{{\bf NOTE}}
		This covers all  sections of \cite{scholze-stix} except for notations and generalities (and also covers the main assertion of \cite{scholze-review}). 
	}
\end{rosetta}

\newpage
\subsection{Final Conclusion}\label{ss:conclusion}
\newcommand{\abc}{$abc$-conjecture}
I have examined the claims of \iut, \cite{scholze-stix} and \cite{scholze-review} in meticulous detail, and provided, in \present, a canonical formulation of the theory without which Mochizuki's work stands incomplete (see \ref{se:why}).

I have also conclusively shown that all the objections voiced in \cite{scholze-stix,scholze-review} regarding Mochizuki's Inter-Universal Teichm\"uller Theory (\iut)  stand dismantled. To put it simply: Table~\ref{tab:rosetta-stone2} shows that every mathematical assertion of \cite{scholze-stix}, \cite{scholze-review} is mathematically incorrect. 

This report also applies to, and invalidates, the discussions of Mochizuki's work on \href{https://www.math.columbia.edu/~woit/wordpress/}{Peter Woit's Blog}, especially but not limited to \cite{woit-summary}, and assertions made elsewhere by some others echoing \cite{scholze-stix}, \cite{scholze-review}.

Hence,   \textit{at this point in time},
the \abc\ stands established by the combination of  \cite{joshi-teich,joshi-untilts,joshi-teich-estimates,joshi-teich-def,joshi-teich-rosetta,joshi-teich-abc-conj} and \cite{mochizuki-iut1,mochizuki-iut2,mochizuki-iut3,mochizuki-iut4}. 

[Recently updated (Feb 2025) versions of all my papers on this topic are available at arxiv.org.]

\Address
\end{document}